%% file: qkzquiver.tex
\documentclass[12pt,a4paper]{amsart}
\usepackage{amssymb}
\usepackage[T1]{fontenc}
\usepackage[a4paper,left=2.25cm,right=2.25cm,top=3cm,bottom=3cm]{geometry}
\usepackage{tikz}
\usetikzlibrary{calc,decorations.pathmorphing,decorations.markings,matrix}
\usepackage{kbordermatrix}
\usepackage{hyperref}
\newcommand{\bra}[1]{\left\langle #1\right|}
\newcommand{\ket}[1]{\left|#1\right\rangle}

\renewcommand\ss{\scriptstyle}
\newcommand\sss{\scriptscriptstyle}
\newcommand\Hom{\operatorname{Hom}}
\newcommand\M{\mathfrak M}
\newcommand\Z{\mathfrak Z}
\newcommand\gl{\mathfrak{gl}}
\newcommand\C{\mathbb C}
\newcommand\Q{\mathbb Q}
\newtheorem{theorem}{Theorem}
\newtheorem{proposition}{Proposition}
\newcommand\B{{\mathrm B}}
\renewcommand\i{{\mathrm i}}
\renewcommand\j{{\mathrm j}}
\allowdisplaybreaks[1]
\long\def\rem#1{}
\author{Paul Zinn-Justin}
\thanks{The author is supported by ERC grant 278124 ``LIC'' and ARC grant DP140102201. He would like to thank A.~Knutson for
discussions and comments in the framework of a parallel collaboration, as well as
V.~Gorbunov, D.~Maulik, A.~Okounkov and N.~Reshetikhin for discussions.
Computerized checks of the results of this paper have been performed
with the help of Macaulay2 \cite{M2}.}
\address{Laboratoire de Physique Th\'eorique et Hautes \'Energies, CNRS UMR 7589 and Universit\'e Pierre et Marie Curie (Paris 6), 4 place Jussieu, 75252 Paris cedex 05, France}
\title[Quiver varieties and qKZ]{Quiver varieties and\\quantum Knizhnik--Zamolodchikov equation}
\input tableau.tex
\begin{document}
\begin{abstract}
We show how equivariant volumes of tensor product quiver varieties of type A are given by matrix elements of
vertex operators of centrally extended doubles of Yangians, and how they satisfy in some cases the rational,
level 1, quantum Knizhnik--Zamolodchikov equation. 
\end{abstract}

\maketitle

\section{Introduction}
It is increasingly clear that a
deep connection exists between generalized cohomology theories and quantum integrable systems, as first pointed out
in \cite{GKV-coh}. The simplest example is that rational quantum integrable systems should be related to (equivariant)
cohomology. Such a connection somewhat accidentally reappeared in \cite{artic33,artic34} in the study
of certain solutions of the quantum Knizhnik--Zamolodchikov equation.
One of the objectives of this short paper is to
revisit this connection in view of recent developments, in particular of the work \cite{MO-qg},
which provides a natural framework for it. The most basic ingredient in the definition of a quantum integrable system
is to provide a solution of the Yang--Baxter equation, the so-called $R$-matrix. And indeed,
\cite{MO-qg} gives a ``recipe'' to compute the $R$-matrix in terms of certain algebro-geometric data.
Conversely, in \cite{artic33,artic34}, starting from certain quantum integrable models, algebraic varieties
were introduced in
such a way that the integrable model performed computations in equivariant cohomology on the latter. These
two approaches are closely related, and in fact, among the equations defining the quantum Knizhnik--Zamolodchikov
equation, one finds the so-called exchange relation which is a direct corollary of the definition 
of the $R$-matrix given in \cite{MO-qg}.

In the present work, we extend the study that was performed
in \cite{artic34,artic56} to arbitrary Nakajima quiver varieties of type A. 
In section 2, we define these quiver varieties and
give an explicit description (only available in type A) which is particularly convenient for explicit computations.
In section 3, we first discuss the exchange relation and the definition of the $R$-matrix. 
We then introduce a key technical tool, which are Vertex Operators for a centrally extended double of a Yangian.
We proceed to show that equivariant volumes of tensor product quiver varieties are equal to matrix elements of
these Vertex Operators, and conclude that in some cases they satisfy the quantum Knizhnik--Zamolodchikov equation.
Since the presentation is rather sparse on details, an example is fully worked out in appendix~\ref{app:ex}.
In particular, this paper is not meant to be fully rigorous mathematically, and some proofs are sketched or 
skipped, trying to focus on the main ideas rather than on technicalities.

\section{Quiver varieties}\label{sec:quiv}
We briefly define here the relevant geometric objects, and refer to the abundant literature for details, 
in particular the work of Nakajima \cite{Nakaj-quiv1,Nakaj-quiv2,Nakaj-quiv3,Nakaj-quiv4} and the lecture notes \cite{Ginz-naka}
and references therein.

\subsection{Definition}
Given two sequences of finite-dimensional vector spaces $V_a$ and $W_a$, with
dimensions $v=(v_a)$ and $w=(w_a)$, $a=1,\ldots,k-1$, we define the corresponding Nakajima quiver variety
of type $\mathrm{A}_{k-1}$
as the quotient of the space of linear maps $(\B_{a,a+1},\B_{a+1,a},\i_a,\j_a)$ in
\[
\bigoplus_{a=1}^{k-2} (\Hom(V_a,V_{a+1})\oplus \Hom(V_{a+1},V_a))
\oplus
\bigoplus_{a=1}^{k-1} (\Hom(W_a,V_a)\oplus \Hom(V_a,W_a))
\]
subject
to the moment map conditions $\mu=0$, where
\[
\mu=(\B_{a,a-1}\B_{a-1,a}-\B_{a,a+1}\B_{a+1,a}+\i_a\j_a)_{a=1,\ldots,k-1}
\]
(where conventionally an index $k$ or $0$ means that the term is absent)
by the natural action of $G_v=\prod_{a=1}^{k-1} GL(V_a)$, namely
\[
(g_a)\in G_v: (\B_{a,a+1},\B_{a+1,a},\i_a,\j_a)\mapsto (g_{a+1} \B_{a,a+1} g_a^{-1},g_a \B_{a+1,a} g_{a+1}^{-1},g_a \i_a, \j_a g_a^{-1})
\]

The algebro-geometric quotient is obtained by considering the subring of $G_v$-invariant functions:
\[
\M_0(v,w)
=\mu^{-1}(0)
/\!/G_v
\]
$\M_0(v,w)$ is an affine singular variety.

Most of the time $\M_0(v,w)$ will be sufficient for our purposes, but at one point we shall need to consider
the GIT quotient
\[
\M(v,w)
=\mu^{-1}(0)
/\!/_\theta\, G_v
\]
where $\theta$ is some appropriate character, see e.g.~\cite{Nakaj-quiv3,Ginz-naka}.
$\M(v,w)$ is smooth quasi-projective.

When there is no risk of confusion, we simply write $\M_0=\M_0(v,w)$ and so on.

\subsection{Combinatorial data}
Associate to $w$ and $v$ the $SL(k)$ weights $\mu=\sum_{a=1}^{k-1} w_a \omega_a$
and $\lambda=\sum_{a=1}^{k-1} w_a \omega_a -\sum_{a=1}^{k-1}v_a\alpha_a$, where the $\omega_a$ are the fundamental
weights of $SL(k)$ and $\alpha_a$ are its simple roots. We shall assume in what follows that
$\lambda$ is a weight of the irreducible representation of $SL(k)$ with highest weight $\mu$.

Also, we shall need to lift $\mu$ and $\lambda$ to $GL(k)$ weights. This is not unique,
but there is an obvious way which is to leave the $k^{\textrm th}$ weight of $\mu$ to be zero.\footnote{The effect of making another choice of lift to $GL(k)$ will be discussed in sect.~\ref{sec:recwheel}.}
The resulting $GL(k)$ weights are conveniently described as arrays of boxes;
as an example,
\begin{alignat*}{3}
w&=(2,1,3) & \quad&\to\quad
\mu=2\tableau{\\\vdotscell\\\vdotscell\\\vdotscell\\}+\tableau{\\\\\vdotscell\\\vdotscell\\}+3\tableau{\\\\\\\vdotscell\\}
&&=
\begin{tikzpicture}[baseline=0]
\node[tableau] {&&&&&\\&&&\\&&\\\vcell\\};
\end{tikzpicture}
\\
v&=(1,0,1) &\quad&\to\quad
\lambda=\begin{tikzpicture}[baseline=0]
\node[tableau] {&&&&&\\&&&\\&&\\\vcell\\};
\draw[bend left=60,->] (tab-1-6.east) ++(0.1,0) to node[right] {$1$} ++(0,-0.5);
\draw[bend left=60,->] (tab-3-3.east) ++(0.1,0) to node[right] {$1$} ++(0,-0.5);
\end{tikzpicture}
&&=
\begin{tikzpicture}[baseline=0]
\node[tableau] {&&&&\\ &&&&\\&\\\\};
\end{tikzpicture}
\end{alignat*}

Moreover, for the sake of simplicity, we shall also assume in what follows that $\lambda$ is dominant\footnote{This hypothesis
can be easily dispensed with by replacing in what follows $\M_0(v,w)$ with the image of $p$.}
The natural morphism $p: \M(v,w) \to \M_0(v,w)$
is then surjective, and a resolution of singularities of $\M_0(v,w)$.

\subsection{Transverse slice of nilpotent orbit closures}\label{sec:mv}
There are several alternative description of $\M_0$ 
in terms of transverse slices of nilpotent orbit closures.
One is particularly convenient
for our purposes: the slice of Mirkovi\'c and Vybornov \cite{MVy-slice}.
\rem{why particularly appropriate: (1) mdeg should be invariant [up to normalization] of embedding. (2) the fusion argument really only relies on the location of the 1's.... (3) only thing I can think of is torus action: maybe the other embedding doesn't make the torus action as obvious??? (4) AH! slodowy slice NOT transverse in $\Z$!!!!}

Denote $N=\sum_{a=1}^{k-1} w_a$, $M=\sum_{a=1}^{k-1} a \,w_a$.
Starting from $\mu$ and $\lambda$, we associate 
two new sequences of integers, $m=(m_1,\ldots,m_N)$ and $\ell=(\ell_1,\ldots,\ell_N)$ respectively,
as follows (see also \cite[sect.~5.1]{MVy-slice}). They are an arbitrary permutation of the sequence of numbers of boxes
in each column of their box diagram as defined in previous section (where for simplicity
of notation we always take $\ell$ to be of length $N$, i.e., possibly pad it with zeroes).
Let us immediately mention that the ordering of $\ell$ is irrelevant; however,
the ordering of $m$ leads to isomorphic but distinct constructions, and it is convenient to retain this freedom.

In the example above, $m$ and $\ell$ are permutations of $(3,3,3,2,1,1)$ and
$(4,3,2,2,2,0)$, respectively.

Note $\sum_{i=1}^N m_i=\sum_{i=1}^N \ell_i=M$.

Then $\M_0$ is isomorphic to $\overline{\mathcal O_\ell} \cap \mathcal T_m$,
where 
\begin{itemize}
\item $\mathcal O_\ell$ is the $GL(M)$-orbit of nilpotent 
$M\times M$ matrices of Jordan type $\ell$, i.e.,
with Jordan blocks of size $\ell_i$,
$i=1,\ldots,N$.

\item $\mathcal T_m=x_m+\mathcal T'_m$,
\begin{itemize}
\item $x_m$ is a nilpotent operator of Jordan type $m$, which we shall choose to be in Jordan form,
i.e., $x_m$ is block diagonal with blocks $\pi_1=[1,m_1]$, $\pi_2=[m_1+1,m_1+m_2]$, $\ldots$, $\pi_N=[M-m_N+1,M]$
and $1$'s right above the diagonal in each diagonal block.
\item $\mathcal T'_m$ is a certain
linear subspace defined as follows: in each block $\pi_i\times\pi_j$
(with the same decomposition of rows and columns
into intervals as above), the only
nonzero entries are on the last row and $\min(m_i,m_j)$ leftmost columns.
\end{itemize}
\end{itemize}

An important property of the intersection of $\overline{\mathcal O}_\ell$ and $\mathcal T_m$ is that
it's {\em transverse}, see \cite{MVy-slice} for details.

\rem{really, it we don't talk about Howe duality, the fact that $\lambda$ and $\mu$ have (at most for $\lambda$;
implied by being in convex envelope of $W\cdot\mu$)
$N$ columns is irrelevant}

\subsection{Tensor product quiver varieties}\label{sec:tpquiv}
Suppose the direct sum $W:=\bigoplus_{a=1}^{k-1} W_a$ 
is given an (ordered) basis $(e_1,\ldots,e_N)$, such that $e_i\in W_{m_i}$ for some 
$m\in \{1,\ldots,k-1\}^N$ with $\#\{i: m_i=a\}=w_a$. Note that the latter condition characterizes precisely
the sequences $(m_1,\ldots,m_N)$ of the previous section.

This provides an isomorphism between $\gl(W)$ (resp.\ $GL(W)$) and $\gl(N)$ (resp.\ $GL(N)$), 
as well as a Borel subgroup, with its unipotent subgroup, Cartan torus and their Lie algebras. 
Explicitly, via the isomorphisms above, we define:
\begin{align*}
\mathfrak n &=\{\text{strict upper triangular matrices}\} \subset \gl(W)
\\
T&=\{\text{invertible diagonal matrices}\} \subset GL(W)
\end{align*}
Note that in fact, $T\subset \prod_{a=1}^{k-1} GL(W_a)$, so it acts in the natural way: $g\in T:(\B_a,\i_a,\j_a)\mapsto (\B_a,\i_a g^{-1},g \j_a)$, and since this action commutes with that of $G_v$, acts on $\M_0$ (or $\M$).

Consider the particular one parameter subgroup $g(t)=\text{diag}(t^{-i},\,i=1,\ldots,N)\in T$, $t\in{\C}^\times$. Then define
\[
\Z_0 = \{ x\in \M_0:\, g(t)x\buildrel t\to0 \over\to 0\}
\]
At one point, we shall also need to consider $\Z:=p^{-1}(\Z_0)$.

Using the alternative description of sect.~\ref{sec:mv}, we can reexpress $\Z_0$ as follows.
First embed $\gl(N)$ inside $\gl(M)$ as matrices that are proportional to the identity inside
each block $\pi_i\times \pi_j$ corresponding to the Jordan decomposition of $x_m$. Then via this embedding,
we simply have
\[
\Z_0 = \M_0 \cap \mathfrak{n}
\]

The affine scheme $\Z_0$ is typically reducible,
and equidimensional of dimension $\frac{1}{2}\dim\M_0$ (this is of course related
to an underlying Lagrangianness, which we do not discuss here).
We denote $\Z_{0,\alpha}$ its irreducible components,
where $\alpha$ runs over some indexing set, which can be chosen as follows.
It consists of the set of tableaux of shape $\lambda$ with the following restrictions:
\begin{itemize}
\item The letter $i$ is used $m_i$ times, $i=1,\ldots,N$.
\item Rows are strictly increasing; columns are weakly increasing (i.e., it is the conjugate of a semi-standard Young tableau).
\end{itemize}
This labelling can be obtained using a generalization of Spaltenstein's algorithm \cite{Spal},
which consists in taking a generic element $X$ of $\Z_{0,\alpha}$ viewed as a $M\times M$ matrix 
and computing the Jordan forms of upper-left submatrices $(X_{ij})_{1\le i,j \le m_1+\cdots+m_h}$,
$h=1,\ldots,N$; the corresponding increasing sequence of partitions forms the tableau $\alpha$.

Similarly, one can introduce $\Z=p^{-1}(\Z_0)$. Its irreducible components $\Z_\alpha$
can be indexed as follows: $\alpha$ is a collection of subsets $\alpha_i$ of $\{1,\ldots,k\}$, $i=1,\ldots,N$,
such that $\# \alpha_i=m_i$ and $\#\{i:a\in\alpha_i\}=\lambda_a$, $a=1,\ldots,k$
(where the $\lambda_a$ are the lengths of rows of $\lambda$ viewed as a box diagram).
Then the map $\varphi$ from $Irr(\mathfrak \Z_0)$ to $Irr(\mathfrak \Z)$ corresponding
to taking $p^{-1}$ is to associate to a tableau $T$ 
the sequence $\alpha_i=\{\text{rows where $i$ is in $T$}\}$.
This map is in general not surjective.

\subsection{Multidegrees}\label{sec:mdeg}
As noted above, the torus $T\cong ({\C}^\times)^N$ acts on $\M_0$. Furthermore, an extra $\mathbb{C}^\times$ action
can be introduced as follows. In the original definition of $\M_0$, it simply scales every linear map:
$\mathbb{C}^\times\ni t:(\B_{a,a+1},\B_{a+1,a},\i_a,\j_a)\mapsto(t \B_{a,a+1},t\B_{a+1,a},t\i_a,t\j_a)$. In the picture of sect.~\ref{sec:mv}, it acts on the slice $\mathcal T_m$ by scaling 
the entry on the last row of each block $\pi_i\times\pi_j$ with $t^{\text{perimeter}-2(\text{column}-1)}$,
where the perimeter is $m_i+m_j$, and the column goes from $1$ to $\min(m_i,m_j)$.

It is not hard to check that $T\times\C^{\times}$ leaves $\Z_0$, and therefore its components $\Z_{0,\alpha}$, invariant.
We want to consider their {\em equivariant volumes}\/ with a particular polynomial normalization. Recall that
from sect.~\ref{sec:mv} we have an embedding, which we denote $\iota$, from $\Z_0=\M_0\cap \mathfrak{n}$ to $\mathcal T_m \cap\mathfrak{n}$. We can then consider pushforward in $(T\times\C^{\times})$-equivariant cohomology:
\[
\Psi_\alpha:=\iota_\ast[\mathfrak Z_{0,\alpha}]
\in H^\ast_{T\times \C^\times}(\mathcal T_m \cap \mathfrak n)
\]
\rem{argh not smooth, so BM homology and cohomology not obviously related by poincare... need to think}
\rem{could we just take equivariant volumes of the $\Z_\alpha$????? should it make any difference? what
about in K-theory?}

Since $\mathcal T_m\cap\mathfrak n$ is equivariantly contractible, its equivariant cohomology is that of a point, i.e.,
a polynomial ring. We write it as
\[
H^\ast_{T\times \C^\times}(\mathcal T_m \cap \mathfrak n)
=
\mathbb Z[z_1,\ldots,z_N,\hbar/2]
\]
where each $z_i$ is the coordinate on $\mathfrak{t}:=\text{Lie}(T)$ corresponding to the basis element $e_i$,
and $\hbar/2$ corresponds to the extra $\C^\times$ action.\footnote{The factor $1/2$ in the generator $\hbar/2$ is there
to match the conventions of integrable system. In fact, one should 
further perform the substitution $\hbar\to -\hbar$ for an exact match; 
but that would obscure the notion of positivity coming from geometry.}
The $\Psi_\alpha$ are sometimes called the {\em multidegrees}\/ of the $\Z_{0,\alpha}$ in $\mathcal T_m \cap \mathfrak n$.
They are homogeneous polynomials of degree, the codimension of $\Z_0$ in $\mathcal T_m \cap \mathfrak n$.
The latter is independent of $m$ (this will be a consequence of transversality, see sect.~\ref{sec:fus}), 
and is given by
$\text{codim}\, \Z_0 = \sum_{a=1}^k \frac{\lambda_a(\lambda_a-1)}{2}$
(where the $\lambda_a$ are the lengths of the rows of the box diagram of $\lambda$).

\subsection{Representation-theoretic interpretation}\label{sec:RT}
Denote by $L_\alpha$ the irreducible representation of $SL(k)$ with highest weight $\alpha$.
It is known that there is a natural isomorphism between the top Borel--Moore homology of $\Z_0$
and $\Hom_{SL(k)}(L_\lambda,L_{\omega_{m_1}}\otimes L_{\omega_{m_2}}\otimes\cdots\otimes L_{\omega_{m_N}})$.
In particular the number of irreducible components of $\Z_0$ is equal to the multiplicity of
$L_\lambda$ in $L_{\omega_{m_1}}\otimes L_{\omega_{m_2}}\otimes\cdots\otimes L_{\omega_{m_N}}$.

We can therefore reinterpret $\Psi:=\iota_\ast$ as an element of $\mathcal H\otimes \mathbb Z[z_1,\ldots,z_N,\hbar/2]$
where
\[
\mathcal H:=
(L^\ast_{\omega_{m_1}}\otimes L^\ast_{\omega_{m_2}}\otimes\cdots\otimes L^\ast_{\omega_{m_N}}
\otimes L_\lambda)^{SL(k)}
\]
More explicitly, $\Psi=\sum_\alpha \Psi_\alpha\, u^\alpha$ where the $u^\alpha$
form the basis dual to the $[\mathfrak Z_{0,\alpha}]$.

\section{Vertex Operators and quantum Knizhnik--Zamolodchikov equation}
In the previous section, given a quiver variety and a basis $(e_1,\ldots,e_N)$ of 
$W=\bigoplus_{a=1}^{k-1} W_a$, 
we have defined geometrically a set of polynomials $\Psi_\alpha$ which
we have grouped into a polynomial-valued vector $\Psi\in\mathcal H\otimes \mathbb Z[z_1,\ldots,z_N,\hbar/2]$. We now study properties of these polynomials.

We shall emphasize the dependence of our construction on the choice of basis $(e_1,\ldots,e_N)$ as follows. We shall denote when required
$\Z_0 := \Z_0^{e_1,\ldots,e_N}$, and similarly for other geometric objects. Note however
that $\Psi$ only depends on the $e_i$ via the combinatorial data of the $m_i$
(such that $e_i\in W_{m_i}$); indeed given two bases with the same $m_i$, there is an obvious
isomorphism from each $W_a$ to itself which is equivariant w.r.t.\ the corresponding two tori,
and so the equivariant classes of the $\Z_{0,\alpha}$ viewed as elements of
$\mathbb Z[z_1,\ldots,z_N,\hbar/2]$ are the same. \rem{literally the same,
not permutation of variables! footnote on that?} We therefore denote
\[
\Psi=\Psi^{m_1,\ldots,m_N},
\quad \mathcal{H}=\mathcal{H}^{m_1,\ldots,m_N}
\]

\subsection{Exchange relation}
The first important relation that $\Psi$ satisfies is the so-called exchange relation.
As we shall see, it is a set of relations involving not a single $\Psi^{m_1,\ldots,m_N}$ 
(unless the $m_i$'s are all equal), but rather all permutations of the sequence of $m_i$'s.
This section is closely related to \cite[chapter 4]{MO-qg} in the sense
that what is called here exchange relation is essentially the {\em definition}\/ of the
$R$-matrix there, with some subtleties, as we shall see in the proof of the following
\begin{proposition}\label{prop:exch}
There exist $\check R_i(z)\in \Hom(
\mathcal{H}^{m_1,\ldots,m_{i},m_{i+1},\ldots,m_N},
\mathcal{H}^{m_1,\ldots,m_{i+1},m_i,\ldots,m_N})\otimes \Q(z,\hbar)$, 
$i=1,\ldots,N-1$,
which satisfy
\begin{align}\label{eq:ybe}
\check R_i(u)
\check R_{i+1}(u+v)
\check R_i(v)
&=
\check R_{i+1}(v)
\check R_i(u+v)
\check R_{i+1}(u)
&&
i=1,\ldots,N-2
\\\label{eq:unit}
\check R_i(u)\check R_i(-u)&=1
&&
i=1,\ldots,N-1
\\\label{eq:comm}
\check R_i(u)
\check R_j(v)
&=
\check R_i(v)
\check R_j(u)
&&
|i-j|>1
\end{align}
such that
\begin{multline}\label{eq:exch}
\Psi^{m_1,\ldots,m_{i+1},m_i,\ldots,m_N}(z_1,\ldots,z_{i+1},z_i,\ldots,z_N)
\\
= \check R_i(z_i-z_{i+1})
\Psi^{m_1,\ldots,m_i,m_{i+1},\ldots,m_N}(z_1,\ldots,z_{i},z_{i+1},\ldots,z_N),
\qquad i=1,\ldots,N-1.
\end{multline}
\end{proposition}
Note that $\check R_i(z)$ is a rational function of both $z$ and $\hbar$, but
it is traditional in integrable systems to emphasize in the notation the dependence on $z$.

We explain here the main steps of the proof:
\begin{itemize}
\item The first basic idea is to compare the various $\Z^{e_1,\ldots,e_N}$
when one permutes the $e_i$ (the appearance of the symmetric group $\mathcal{S}_N$
is related to the ``gauge group'' $\prod_a GL(W_a)$ being of type A). 
The classes of the irreducible components $\Z_\alpha^{e_1,\ldots,e_N}$ 
are known to form a basis of the (appropriately localized)
equivariant cohomology ring of $\M$. Since permutations of the $e_i$'s do not affect
the torus $T$, we have at our disposal $N!$ different bases of $H^{\ast,loc}_{T\times\C^\ast}(\M)$
as free modules over $H^{\ast,loc}_{T\times\C^\ast}(\cdot)$.
We may as well consider elementary transpositions only, and then write formally
\[
[\mathfrak Z^{\ldots e_{i+1},e_{i}\ldots}_\alpha]=\sum_\beta
(\check {\mathfrak R}_i)_{\alpha}^{\beta} [\mathfrak Z^{\ldots e_{i},e_{i+1}\ldots}_\beta]
\]
where the $\check {\mathfrak R}_i)_{\alpha}^{\beta}$ are matrices
in $H^{\ast,loc}_{T\times\C^\ast}(\cdot)$.
This analysis is not entirely satisfactory, and
we mention briefly how it can be improved (this will not be necessary for the statement
of any of the theorems in this paper, but is definitely needed to prove them).

One can perform on the basis of the $[\Z_\alpha]$ any further linear transformation 
that is compatible with the $\mathcal{S}_N$ action. A good choice is to switch to
the ``stable'' basis, which we shall not define here (see \cite{MO-qg} and app.~\ref{app:ex}),
and which corresponds to a change from the $[\Z_\alpha]$ to
the say $[\Z'_\alpha]=\sum_\beta c_\alpha^\beta [\Z_\beta]$. 
The coefficients $c_\alpha^\beta$ of the change
of the basis are upper triangular (w.r.t.\ to a natural order on the indexing set). 
and $\mathbb Z$-valued; and the $R$-matrices
in this new basis, defined by
$[\mathfrak Z'{}^{\ldots e_{i+1},e_{i}\ldots}_\alpha]=\sum_\beta
(\check {\mathfrak R'}_i)_{\alpha}^{\beta} [\mathfrak Z'{}^{\ldots e_{i},e_{i+1}\ldots}_\beta]$,
have various nice properties which are explained in \cite{MO-qg}. In particular, they have
a ``locality'' property w.r.t.\ tensor product. One practical consequence is that
the matrix entries of $\check{\mathfrak R}'_i$ (and therefore of $\check{\mathfrak R}_i$),
as elements of $\Q(z_1,\ldots,z_N,\hbar/2)$, are functions of $z_i-z_{i+1}$ and $\hbar$ only.

\item Next we want to pushforward from $\M$ to $\M_0$ using $p_\ast$.
Since $\Z$ and $\Z_0$ are equidimensional
of same dimension, given a component $\Z_\alpha$ of $\Z$, if one applies $p$
one of two things may happen: either it is sent to something lower-dimensional, 
or it is sent to an irreducible component of $\Z_0$, which is nothing but $\Z_{0,\alpha}$ 
if we identify a tableau $\alpha$ with its image under $\varphi$ (see sect.~\ref{sec:tpquiv});
in which case the map turns out to be generically one-to-one.
This means that in equivariant cohomology, $[\Z_\alpha]$ is sent to either zero
or $[\Z_{0,\alpha}]$. In the end, we get the exact same equality
\begin{equation}\label{eq:preexch}
[\mathfrak Z^{\ldots e_{i+1},e_{i}\ldots}_{0,\alpha}]=\sum_\beta
(\check {\mathfrak R}_i)_{\alpha}^{\beta} [\mathfrak Z^{\ldots e_{i},e_{i+1}\ldots}_{0,\beta}]
\end{equation}
where the range of indices has been appropriately restricted.

\item Finally we want to compute equivariant volumes, with a prescribed normalization
dictated by the embedding space. Two complications need to be addressed at this stage:
\begin{itemize}
\item In order to turn equality \eqref{eq:preexch} into an equality of polynomials,
one has to parameterize the torus. Here we break the symmetry between
$\ldots e_i,e_{i+1}\ldots$ and 
$\ldots e_{i+1},e_{i}\ldots$ by stating that $z_i$ is the coordinate corresponding to $e_i$.
If we then apply $\iota_\ast$ to both sides of \eqref{eq:preexch}, the r.h.s.\ behaves
as expected: $\check {\mathfrak R}_i$ becomes a
$\Q(z_1,\ldots,z_N,\hbar/2)$-valued matrix, and is fact as already mentioned is
$\Q(z_i-z_{i+1},\hbar)$-valued 
(one can be even more precise since the localization only allows
denominators of the form of a product of $a\hbar+z_i-z_{i+1}$, $a\in\mathbb Z$);
and $[\mathfrak Z^{\ldots e_{i},e_{i+1}\ldots}_{0,\beta}]$ becomes $\Psi_\beta^{...m_i,m_{i+1}\ldots}$.
However the l.h.s.\ is not exactly $\Psi_\alpha^{...m_{i+1},m_{i}\ldots}$ as one might naively
expect because of the mismatch in the parameterization of the torus: one needs
to exchange variables $z_i$ and $z_{i+1}$ to take it into account.

\item
The embedding spaces for $\Z_0^{\ldots e_{i},e_{i+1}\ldots}$
and $\Z_0^{\ldots e_{i+1},e_{i}\ldots}$ are not the same: this
means when we pass from \eqref{eq:preexch} to equivariant volumes, one needs
to introduce a corrective factor corresponding to the ratio of equivariant volumes of
$\mathcal T_{\ldots m_i,m_{i+1}\ldots} \cap {\mathfrak n}$
and
$\mathcal T_{\ldots m_{i+1},m_i\ldots} \cap {\mathfrak n}$.
The latter can in fact be computed explicitly and is
\[
\prod_{a=1}^{\min(m_i,m_{i+1})} \frac{a\hbar - z_i + z_{i+1}}{a\hbar - z_{i+1} + z_{i}}
\]
In order to accommodate for this corrective factor,
we introduce a new normalization of the $R$-matrix:
\begin{equation}\label{eq:newnorm}
\check R_i(z):=\prod_{a=1}^{\min(m_i,m_{i+1})} \frac{a\hbar-z}{a\hbar+z}\ \check{\mathfrak R}_i(z)
\end{equation}
(this is a rather minor issue but is important to compare explicitly the results
of \cite{artic34,artic56} and of app.~\ref{app:ex} to those of \cite{MO-qg}).
\end{itemize}
Putting everything together, we find that Eq.~\eqref{eq:exch} is precisely satisfied.

\item Successive changes of basis should lead to
the relations of the symmetric group ${\mathcal S}_N$
for the operators $\check{\mathfrak R}_i$.
However, as explained in the previous paragraph,
one must take into account the possible reparameterization
of the torus due to permutation of the $e_i$. 
In other words, after a first change of basis leading to $\check{\mathfrak R}_i$,
the next change of basis corresponding to elementary transposition $(j,j+1)$
will not be given by $\check{\mathfrak R}_j$ but rather
by $\tau_i\mathfrak{R}_j$, where $\tau_i$ is the operator on $\Q(z_1,\ldots,z_N,\hbar/2)$
that permutes $z_i$ and $z_{i+1}$.
More generally, viewing the $(\check{\mathfrak R}_i)_\alpha^\beta$ as operators of multiplication
on $\Q(z_1,\ldots,z_N,\hbar/2)$, the precise statement is that the $\tau_i \check{\mathfrak R}_i$
satisfy the symmetric group relations.
Rewriting these explicitly by pulling out the $\tau_i$ (and using the fact that the $\tau_i$ 
alone satisfy the same relations) results in 
Eqs.~(\ref{eq:ybe}--\ref{eq:comm}) first for the $\check{\mathfrak R}_i$, then, noting that
transformation \eqref{eq:newnorm} preserves them, for the $\check{R}_i$.
\end{itemize}

Given the polynomials $\Psi$, it is not hard to see that \eqref{eq:exch} defines uniquely the $\check R_i$
for $N\ge3$. So \eqref{eq:exch} can in principle be considered as a (somewhat clumsy) definition of the $R$-matrix,
and it is in that sense that it is equivalent (for $N\ge3$) to the definition of \cite{MO-qg}.

\newcommand\DY{\widehat{DY(\mathfrak{sl}(k))}}
\subsection{Correlators of Yangian Vertex Operators}
The exchange relation \eqref{eq:exch} does not characterize the polynomial-valued vectors
$\Psi$ uniquely. We now give a formula for them which does (up to normalization).

We start from the centrally extended double of the Yangian of $\mathfrak{sl}(k)$,
which we denote $\DY$. It was introduced in
\cite{Kh-Yd} and the case $k=2$ was studied further in \cite{KhLP-Yd}. One should think of it
as a ``rational'' analogue (in the sense of the classification of the solutions of the 
Yang--Baxter equation) of the ``trigonometric'' quantized affine algebras, and its
representation theory is expected to run parallel to that of the latter, though some facts
that we present now have not been derived rigorously in the literature, as far as this author
knows.

The value of the central element in an irreducible representation of $\DY$ is called
the {\em level}.
We are here interested in two types of modules of $\DY$:
\begin{itemize}
\item Evaluation modules $L_\nu(z)$, where $z$ is a formal parameter,
which are level 0 and isomorphic to $L_\nu \otimes \C[z,z^{-1}]$ as $U(\mathfrak{sl}(k))$-modules.
\item Highest weight modules $\mathcal V_{1,\lambda}$, where we separated the
highest weight of $\widehat{\mathfrak{sl}(k)}$ into the level, which for us will always be $1$,
and the highest weight $\lambda$ of 
the finite-dimensional algebra $\mathfrak{sl}(k)$; in fact in level $1$, 
$\lambda$ can only be zero (case of the so-called basic
representation) or a fundamental weight.
\end{itemize}

Next we define type \textrm{I}  (dual) vertex operators (VOs) $\Phi_{\omega_a}(z)$, $a=1,\ldots,k-1$, i.e.,
intertwiners
\[
\Phi^{1,\lambda,\lambda'}_{\omega_a}(z): \mathcal V_{1,\lambda} \otimes L_{\omega_a}(z) \to \mathcal V_{1,\lambda'}
\]
They are unique up to multiplication by a scalar, and the exact choice of normalization is irrelevant for our purposes.
Furthermore, they are ``perfect'' VOs in the sense that $\lambda'$ is entirely determined by $\lambda$ and $\omega_i$. 
Explicitly,
if $\lambda=\omega_b$, $b=0,\ldots,k-1$ (with $\omega_0:=0$), then
$\lambda'=\omega_{b+a \pmod k}$. Because of this fact one can safely suppress the superscripts and write
simply $\Phi_{\omega_a}(z)$.

An important ingredient is the bosonic realization of the level $1$ representations. This construction is performed
in \cite{Kh-Yd,KhLP-Yd} for $k=2$ and can be extended straightforwardly to any $k$.
It consists of a $\DY$-module $\mathcal B$ which is a bosonic Fock space, that is the direct sum of 
copies of the highest weight module of $\widehat{\mathfrak{gl}(1)^{k-1}}$, the Heisenberg algebra sitting inside $\DY$. 
The highest weight vectors are denoted $\ket{\lambda}$,
where $\lambda$ runs over the weight lattice of $\mathfrak{sl}(k)$. In particular if $\lambda$ is zero or
a fundamental weight, $\ket{\lambda}$ is a highest weight vector of the whole of $\DY$, with associated
submodule $\mathcal V_{1,\lambda}$.

The dual module $\mathcal B^\ast$ is similarly a direct sum of lowest weight modules, and we denote the
lowest weight vectors $\bra{\lambda}$.

One more notation is needed: a highest weight vector $\ket{\lambda}$ of $\widehat{\mathfrak{gl}(1)^{k-1}}$, 
where $\lambda$ is dominant,
generates a submodule of the finite-dimensional algebra $\mathfrak{sl}(k)$ which is isomorphic to $L_\lambda$.
This defines a unique (up to normalization) vector of $L_\lambda^\ast \otimes \mathcal B$ which we denote
$\ket{L_\lambda^\ast}$. Below we need $\ket{L_\lambda}$ rather than $\ket{L_\lambda^\ast}$;
note that $L_\lambda^\ast \cong L_{\lambda^\ast}$ where by definition $\lambda^\ast:=-w_0 \lambda$, $w_0$ longest permutation
in $\mathcal S_k$.

We can then state:
\begin{theorem} \label{thm:main}
The following formula holds:
\begin{equation}\label{eq:main}
\Psi(z_1,\ldots,z_N)=
\frac{1}{\kappa(z_1,\ldots,z_N)}
\bra{0}
\Phi_{\omega_{m_1}}(z_1)\ldots \Phi_{\omega_{m_N}}(z_N)
\ket{L_\lambda}
\end{equation}
where $\kappa(z_1,\ldots,z_N)$ is a meromorphic scalar function.
\end{theorem}
Note that the product of VOs in this expression is the product as linear operators on $\mathcal B$.
All together, the r.h.s.\ of \eqref{eq:main} is an element of 
$L^\ast_{\omega_{m_1}}\otimes L^\ast_{\omega_{m_2}}\otimes\cdots\otimes L^\ast_{\omega_{m_N}}
\otimes L_\lambda$, as should be (and in fact the intertwining property of the VOs for $\mathfrak{sl}(k)$
ensures that it lies more precisely in the $\mathfrak{sl}(k)$-invariant subspace, just like the l.h.s.,
cf sect.~\ref{sec:RT}).

This can be considered the main result of this paper.
In the next paragraphs we provide a sketch of proof of it.
It proceeds in two steps: first we prove it in the case $\mu=N\omega_1$; then we extend it to the general case by
using the fusion procedure. Note that contrary to the exchange relation, we do not know of a geometric proof.

\subsubsection{Case of $\mu=N\omega_1$: orbital varieties}
If $\mu=N\omega_1$, $\Z_0=\overline{\mathcal O_\ell}\cap \mathfrak n$ is known as the orbital scheme associated to the nilpotent orbit closure $\M_0=\overline{\mathcal O_\ell}$,
and its irreducible components $\Z_{0,\alpha}$ 
are called {\em orbital varieties}\/ (and $\M$ is the cotangent bundle of an appropriate
partial flag variety). These were already studied in a similar context in
\cite{artic34,artic56}, although the formula of Thm.~\ref{thm:main} in terms of Vertex Operators was not given there.

One of the main results of \cite{artic56} is that $\Psi$ satisfies the exchange relation \eqref{eq:exch},
where the $R$-matrix was identified explicitly; it is given
by
\begin{equation}\label{eq:Rone}
\check R_i(z)=\frac{\hbar- z P_{i,i+1}}{\hbar+z}
\end{equation}
where $P_{i,i+1}$ permutes factors $i$ and $i+1$ of the tensor product
(there is an overall minus sign mismatch compared to \cite{artic56}, which is due to different conventions).
\rem{really? thought it was $U_{-1}(g)$ R-matrix? check signs}

Now the level 1 VOs $\Phi(z):=\Phi_{\omega_1}(z)$ satisfy an exchange relation of their own:
\begin{equation}\label{eq:exchVO}
\check R^\star(z_1-z_2)\Phi(z_1)\Phi(z_2)=\Phi(z_2)\Phi(z_1)
\end{equation}
but with a different normalization of the $R$-matrix: \cite{KhLP-Yd}
\[
\check R^\star(z)=\frac{f(z)}{f(-z)} \check R(z),
\qquad
f(z)=\frac{\Gamma(1-\frac{z}{k\hbar})}{\Gamma(1-\frac{z-\hbar}{k\hbar})}
\]
This implies that if one sets
\begin{equation}\label{eq:kap}
\kappa(z_1,\ldots,z_N)=\prod_{1\le i<j\le N} f(z_j-z_i)
\end{equation}
and defines as in Thm.~\ref{thm:main}
\[
\Xi(z_1,\ldots,z_N)=\frac{1}{\kappa(z_1,\ldots,z_N)}\bra{0}\Phi(z_1)\ldots\Phi(z_N)\ket{L_\lambda}
\]
then $\Xi$ satisfies the exchange relation \eqref{eq:exch} (with the same normalization
of the $R$-matrix as $\Psi$).

Furthermore, the explicit form \eqref{eq:Rone} (identity plus permutation of indices)
of the $R$-matrix shows that given one entry of $\Psi$ in the standard basis
of a weight space of $(\C^k)^{\otimes N}$
(which is identified in the present context with the stable basis), 
all other entries are determined uniquely by repeated application
of the exchange relation \eqref{eq:exch}. 
Explicitly, writing these entries as $\Psi'_\alpha=\sum_\beta c^\beta_\alpha \Psi_\alpha$ and
similarly for $\Xi$,
if we can show that $\Psi'_\alpha=\Xi'_\alpha$ for some $\alpha$, then $\Psi=\Xi$.

On the geometric side, consider the (standard) tableau $\alpha^{(0)}$ of shape $\lambda$
obtained by filling boxes with $1,\ldots,N$ row by row, e.g.,
\tableau{1&2&3\\4&5\\6\\}. Then define the following linear space
\[
\Z_{0,\alpha^{(0)}}=\{U\ N\times N\text{ upper triangular: }
U_{ij}=0\text{ if $i$ and $j$ are on the same row of $\alpha^{(0)}$}
\}
\]
It is elementary to show that $\Z_{0,\alpha^{(0)}}$ as defined above is an irreducible component of $\Z_0$
(it sits in it and has maximal dimension),
and indeed corresponds to the particular tableau $\alpha^{(0)}$. 
We have immediately
\[
\Psi_{\alpha^{(0)}}
=\prod_{\substack{1\le i<j\le N\\\text{$i$ and $j$ on same row of $\alpha^{(0)}$}}} (\hbar+z_i-z_j)
\]
The tableau $\alpha^{(0)}$ is special in that it is
the smallest w.r.t.\ the natural order on tableaux for which the change of basis to the stable basis is
triangular; that means explicitly that
$
\Psi'_{\alpha^{(0)}} = \Psi_{\alpha^{(0)}}
$
where as usual we identify $\alpha^{(0)}$ with its image under $\varphi$,
namely
$
({\underbrace{1,\ldots,1}_{\lambda_1},\ldots,\underbrace{k,\ldots,k}_{\lambda_k}})
$.

On the VO side, one can perform the calculation of the expectation value along the same lines as 
\cite{KhLP-Yd2} (which only treats the $k=2$ case). Alternatively,
one can start from the corresponding calculations for the quantized affine algebra
$U_q(\widehat{\mathfrak{sl}(k)})$, as in \cite[sect.~9]{Nakay-slNqKZ},
and take the limit $\epsilon\to0$ with
$q=e^{-\epsilon\hbar/2}$, $w=e^{\epsilon z}$ (where $w$ stands here for any multiplicative spectral parameter).

We find that if $\alpha^{(0)}$ is any weakly increasing sequence, i.e.,
is of the form above:
$\alpha^{(0)}=
({\underbrace{1,\ldots,1}_{\lambda_1},\ldots,\underbrace{k,\ldots,k}_{\lambda_k}})
$,
then the following identity holds:
\[
\bra{0}\Phi(z_1)_{\alpha^{(0)}_1}\ldots
\Phi(z_N)_{\alpha^{(0)}_N}\ket{\lambda^\ast}
=
\prod_{1\le i<j\le N} f(z_i-z_j)
\prod_{\substack{1\le i<j\le N\\\alpha^{(0)}_i=\alpha^{(0)}_j}} (\hbar+z_i-z_j)
\]
(where $\Phi(z)_a$, $a=1,\ldots,k$, denotes the entry of $\Phi(z)$ in the standard basis of $\C^k$).

Combining the various formulae above, we conclude that $\Psi'_{\alpha^{(0)}}=\Xi'_{\alpha^{(0)}}$, and therefore $\Psi=\Xi$.

\subsubsection{Fusion procedure and transverse slice}\label{sec:fus}
We first explain the fusion procedure and how it allows us to reduce the computation
of the r.h.s.\ of \eqref{eq:main} to the case of the previous section, with $\mu=M\omega_1$.

Given a fundamental weight $\omega_a$ of $SL(k)$, $a=1,\ldots,k-1$, there is a
(unique up to normalization) embedding $p_a$ of $L_{\omega_a}$ into
$L_{\omega_1}^{\otimes a}$ as $U(\mathfrak{sl}(k))$-modules. 
Furthermore, the space $L_{\omega_1}(z-\frac{a-1}{2}\hbar)\otimes
L_{\omega_1}(z-\frac{a-3}{2}\hbar)\otimes
\cdots
\otimes
L_{\omega_1}(z+\frac{a-1}{2}\hbar)$ is reducible for the whole $\DY$-action,
and the image of $p_a$ forms an invariant subspace which is isomorphic to $L_{\omega_a}(z)$.
(noting again that the usual integrable systems sign convention corresponds to $\hbar\to -\hbar$).

Consider then the operator
\begin{equation}\label{eq:fusion}
\Phi(z-\frac{a-1}{2}\hbar)\Phi(z-\frac{a-3}{2}\hbar)\ldots\Phi(z+\frac{a-1}{2}\hbar)
\,
p_a
\end{equation}
It is by definition an intertwiner,
and we conclude that it is equal to $\Phi_{\omega_a}(z)$ up to normalization,
by uniqueness of this intertwiner.
Therefore the r.h.s.\ of \eqref{eq:main} is nothing but
\begin{multline}\label{eq:rhsfus}
\frac{1}{\kappa'(z_1,\ldots,z_N)}
\bra{0}
\Phi(z_1-\frac{m_1-1}{2}\hbar)\ldots\Phi(z_1+\frac{m_1-1}{2}\hbar)\,p_1
\\
\ldots
\Phi(z_N-\frac{m_N-1}{2}\hbar)\ldots\Phi(z_N+\frac{m_N-1}{2}\hbar)\,p_N
\ket{L_\lambda}
\end{multline}
where we absorbed normalization issues in the scalar function $\kappa'$,
i.e., $\kappa'(z_1,\ldots,z_N)=\kappa(z_1,\ldots,z_N) \prod_{i=1}^N c_{m_i}(z_i)$
(we shall not need to compute the $c_a(z)$ explicitly). 
So the r.h.s.\ now looks like its special case where $\mu=M\omega_1$.

Now let us consider the l.h.s.\ of \eqref{eq:main}. We shall finally
use the construction of sect.~\ref{sec:mv}. Recall that $\M_0 = \overline{\mathcal O_\ell}
\cap \mathcal T_m$, where we reinterpret in the present context 
$\M_0^{(1)}:=\overline{\mathcal O_\ell}$ as the singular Nakajima variety associated to
$\mu=M\omega_1$ and same $\lambda$. Similarly, $\Z_0 = \Z_0^{(1)}\cap \mathcal T_m$, 
where $\Z_0^{(1)}:=\overline{\mathcal O_\ell}\cap\mathfrak{n}$ is the
tensor product quiver variety for $\mu=M\omega_1$, i.e., $m=\underbrace{(1,\ldots,1)}_M$.

A nice property of the Mirkovic--Vybornov slice is that it is still transverse
on $\Z_0^{(1)}$. More specifically, given an irreducible component $\Z_{0,\alpha}^{(1)}$,
the intersection is either empty, or transverse and irreducible 
(irreducibility is in fact also a corollary of the argument that follows).

Now let us compute multidegrees. The $\Z_{0,\alpha}^{(1)}$ are invariant under a torus
$T^{(1)}\cong (\C^\times)^M \times \C^\times$, but the slice with $\mathcal T_m$ only
preserves a subtorus $T\subset T^{(1)}$, with $T\cong (\C^\times)^M \times \C^\times$.
The subtorus $(\C^\times)^M$ is easy to understand: in accordance with sect.~\ref{sec:tpquiv},
it corresponds to conjugation w.r.t.\ invertible
diagonal matrices which are proportional to the
identity in each Jordan block of $x_m$. 
The extra $\C^\times$ is more subtle;
an explicit formula was given at the beginning of sect.~\ref{sec:mdeg}.
We recover it now as follows.

Any element of $\mathcal T_m$, just like
$x_m$ itself, has a $1$ right above the diagonal in each Jordan block. In order for $T$ to
preserve $\mathcal T_m$, the weight of those entries w.r.t.\ $T$ must be zero.
Now the weights of $\mathfrak{gl}(M)$ w.r.t.\ $T^{(1)}$ are simply
$\hbar+z^{(1)}_i-z^{(1)}_j$, $i,j=1,\ldots,M$,
(where we added the superscripts $(1)$ for convenience), each one corresponding to the entry at $(i,j)$.
We conclude that $\mathfrak t^\ast$, viewed as a quotient of 
$\mathfrak t^{(1)}{}^\ast$, is generated by $z^{(1)}_1,\ldots,z^{(1)}_M,\hbar$ with relations
\[
\hbar+z^{(1)}_i-z^{(1)}_{i+1}=0,\qquad i,i+1\in \pi_h,\quad h=1,\ldots,N
\]
Note that the $\C^\times$ itself viewed as a subgroup of $T$ is not canonically 
defined and that corresponds to the freedom in solving the equations above as
$z^{(1)}_i=z_h+cst+ i\hbar$, $i\in\pi_h$, 
to shift $z_h$ by a constant. We choose, as in sect.~\ref{sec:mdeg},
to fix the constant by asking that the average of $z^{(1)}_i$ in each $\pi_h$ be $z_h$.

We finally conclude that
\begin{equation}\label{eq:fusmdeg}
[\Z^{(1)}_{0,\alpha}]_T = \Psi^{(1)}_\alpha
(z_1-\frac{m_1-1}{2}\hbar,\ldots,z_1+\frac{m_1-1}{2}\hbar,
\ldots,
z_N-\frac{m_N-1}{2}\hbar,\ldots,z_N+\frac{m_N-1}{2}\hbar)
\end{equation}
where the l.h.s.\ is the multidegree w.r.t.\ $T$, whereas $\Psi^{(1)}_\alpha(z_1,\ldots,z_M)$
is by definition the multidegree w.r.t.\ $T^{(1)}$.

Now transversality of the intersection means that the multidegree of $\Z^{(1)}_{\alpha,0}$
in $\mathfrak n$ is equal to the multidegree of $\Z^{(1)}_{0,\alpha}\cap \mathcal T_m$
in $\mathfrak n\cap\mathcal T_m$. We conclude that the l.h.s.\ of \eqref{eq:fusmdeg}
is nothing but
$\sum_\beta x_{\alpha,\beta} \Psi_\beta$, where $x_{\alpha,\beta}$ is the multiplicity
of $\Z_{0,\beta}$ in $\Z^{(1)}_{0,\alpha}\cap\mathcal T_m$.

A more detailed analysis, which we shall skip here, 
shows that the stable basis can be identified with the standard basis of (the weight
space of) $L_{\omega_{m_1}}\otimes\cdots\otimes L_{\omega_{m_N}}$
(note that this is essentially forced on us by the locality of the basis w.r.t.\ 
tensor product and the fact that it respects weight space decomposition).
This effectively fixes the matrix $x_{\alpha,\beta}$, and implies that
\begin{multline}
\Psi(z_1,\ldots,z_N)
=\Psi^{(1)}
(z_1-\frac{m_1-1}{2}\hbar,\ldots,z_1+\frac{m_1-1}{2}\hbar,
\\
\ldots,
z_N-\frac{m_N-1}{2}\hbar,\ldots,z_N+\frac{m_N-1}{2}\hbar)
(p_1\otimes\cdots\otimes p_N)
\end{multline}
\rem{would be nice to check explicitly this basis issue}

Comparing this with \eqref{eq:rhsfus}, and using \eqref{eq:main} for $\mu=M\omega_1$,
we fix $\kappa'(z_1,\ldots,z_N):=\kappa(z_1-\frac{m_1-1}{2}\hbar,\ldots,z_1+\frac{m_1-1}{2}\hbar,
\ldots,
z_N-\frac{m_N-1}{2}\hbar,\ldots,z_N+\frac{m_N-1}{2}\hbar)$ with $\kappa$ given by \eqref{eq:kap},
and obtain \eqref{eq:main} in the general case.

As a final remark, the vertex operators $\Phi_{\omega_i}(z)$ satisfy an exchange relation as well:
\[
\check R^\star_{\omega_a,\omega_b}(z_1-z_2)\Phi_{\omega_a}(z_1)\Phi_{\omega_b}(z_2)=\Phi_{\omega_b}(z_2)\Phi_{\omega_a}(z_1)
\]
where from the fusion construction (cf \eqref{eq:exchVO} and \eqref{eq:fusion}),
$\check R^\star_{\omega_a,\omega_b}(z)$ is up to normalization a product of $ab$ $\check R^\star$ matrices.
This identifies the geometrically defined $R$-matrices of Prop.~\ref{prop:exch} with
the $R$-matrices of $\DY$ in the tensor product of two evaluation modules corresponding to fundamental representations,
as expected (as usual, up to a scalar factor).

\subsubsection{Recurrence relations and wheel conditions}\label{sec:recwheel}
As an aside, we can apply the same strategy as in the previous section
to investigate the following question: what happens when in the construction of $\M_0$
as a transverse slice, as in sect.~\ref{sec:mv}, we use another lift of $\mu$
from $SL(k)$ to $GL(k)$, adding extra columns of length $k$? That is, we write
$\tilde m=(m_1,\ldots,m_{p-1},k,m_p,\ldots,m_N)$ and similarly $\tilde\ell=(k,\ell_1,\ldots,\ell_N)$ 
(recall that the ordering of $\tilde\ell$ does not matter).

It is not hard to see that in $\overline{\mathcal O_{\tilde\ell}}\cap \mathcal T_{\tilde m}$.
the $k$ rows and columns of the $M\times M$ matrix corresponding to the block of size $k$ in $x_{\tilde m}$ are
completely fixed, i.e., the diagonal block is of Jordan form and the rest is zero. This provides
a direct isomorphism
$
\overline{\mathcal O_{\tilde\ell}}\cap \mathcal T_{\tilde m}
\cong
\overline{\mathcal O_{\ell}}\cap \mathcal T_{m}=\M_0
$, and similarly for the tensor product variety,
$
\overline{\mathcal O_{\tilde\ell}}\cap \mathcal T_{\tilde m} \cap \tilde{\mathfrak n}
\cong
\overline{\mathcal O_{\ell}}\cap \mathcal T_{m}\cap \mathfrak n=\Z_0
$. Such isomorphisms of varieties in different sizes lead to {\em recurrence relations}
for $\Psi$. As a first step, denote $\Psi^{\tilde m}_{\tilde\alpha}$ the multidegrees of the irreducible components
of $\overline{\mathcal O_{\tilde\ell}}\cap \mathcal T_{\tilde m} \cap \tilde{\mathfrak n}$,
which are naturally indexed by tableaux of $\tilde\lambda$ (the box diagram associated to $\tilde\ell$).
Then from the above,
\begin{multline*}
\Psi^{\tilde m}_{\tilde\alpha}(z_1,\ldots,z_{p-1},\zeta,z_p,\ldots,z_N)
=
\\
\prod_{i=1}^{p-1} \prod_{a=0}^{m_i-1}\left(\Big(\frac{m_i+k}{2}-a\Big)\hbar+z_i-\zeta\right)
\prod_{i=p}^N \prod_{a=0}^{m_i-1}\left(\Big(\frac{m_i+k}{2}-a\Big)\hbar+\zeta-z_i\right)
\Psi_\alpha^m(z_1,\ldots,z_N)
\end{multline*}
where the nilpotent orbit is $\tilde\ell$ in the l.h.s., $\ell$ in the r.h.s.,
and the linear factors take care of the difference of embedding spaces. 
The tableau $\alpha$ is obtained from $\tilde\alpha$ by removing the boxes labelled $p$ and then renumbering
$i\mapsto i-1$ for $i>p$.

One can then view the slice with $\tilde m$ as a further slice of the one with the sequence $(m_1,\ldots,m_{p-1},
n_1,\ldots,n_r,m_p,\ldots,m_N)$ with $n_1+\cdots+n_r=k$,
and where now the $n_i<k$. Comparing the multidegrees leads to an actual recurrence relation
between $\Psi$'s as defined in sect.~\ref{sec:quiv}:
\begin{multline}\label{eq:recur}
\Psi_{\tilde\alpha}^{m_1,\ldots,n_1,\ldots,n_r,\ldots,m_N}
(z_1,\ldots,z_{p-1},\zeta_1,\ldots,\zeta_r,
z_p,\ldots,z_N)\big\vert_{\zeta_{i+1}-\zeta_i=\frac{n_i+n_{i+1}}{2}\hbar}
\\
=\begin{cases}
0&\text{ if }\exists\, i<j\in\{p,p+1,\ldots,p+r-1\},\ \text{row}(i)\ge\text{row}(j) \text{ in }\tilde\alpha
\\
\rlap{$\prod_{i=1}^{p-1} \prod_{a=0}^{m_i-1}\left(\Big(\frac{m_i+n_1}{2}-a\Big)\hbar+z_i-\zeta_1\right)
\prod_{i=p}^N \prod_{a=0}^{m_i-1}\left(\Big(\frac{m_i+n_r}{2}-a\Big)\hbar+\zeta_r-z_i\right)$}
\\
\qquad\Psi_\alpha^{m_1,\ldots,m_N}(z_1,\ldots,z_N)&\text{ otherwise}
\end{cases}
\\[-5.5mm]
n_1+\cdots+n_r=k
\end{multline}
where once again $\alpha$ is obtained from $\tilde\alpha$ by removing boxes labelled $p,\ldots,p+r-1$ and then shifting
$i\mapsto i-r$ for $i\ge p+r$.
(in the case $m_i=n_i=1$, this recurrence relation
was already mentioned in \cite{artic34} but without any geometric proof.)

The reasoning can be pushed one step further. What about if we have an entry (strictly) greater than $k$
in $m$, keeping $\ell$ of the usual form ($\ell_i\le k$ for all $i$)?
In that case we immediately find $\overline{\mathcal O}_\ell \cap \mathcal T_m = \text{\O}$,
and with the same argument, obtain
\begin{equation}\label{eq:wheel}
\Psi^{\ldots,n_1,\ldots,n_2,\ldots,n_r,\ldots}
(\ldots,\zeta_1,\ldots,\zeta_2,\ldots,\zeta_r,\ldots)
\big\vert_{\zeta_{i+1}-\zeta_i=\frac{n_i+n_{i+1}}{2}\hbar}
=0,
\qquad n_1+\cdots+n_r>k
\end{equation}
Here, in contradistinction with \eqref{eq:recur}, the location of the $\zeta_i$ is arbitrary, as long as
they are ordered. Equations \eqref{eq:wheel} are known as {\em wheel conditions}.

\subsection{Rational quantum Knizhnik--Zamolodchikov equation}
We finally discuss an important special case for $\Psi$, which is when $\ket{\lambda}$
is a highest weight vector of $\DY$. As already mentioned, this implies that $\lambda$ is either zero
or a fundamental weight of $SL(k)$; equivalently, it means that
the box diagram of $\lambda$
is either a rectangle of height $k$, or a ``quasi-rectangle'' (a rectangle with one extra
partial column). Note that for a given $w$, there is a unique $v$ such that $\lambda$
satisfies these conditions.

In that case, it is expected that $\Psi(z_1,\ldots,z_N)$ satisfies the rational
quantum Knizhnik--Zamolodchikov equation (also known as difference KZ equation)
\cite{TV1}, which can be considered as the rational limit of the (trigonometric)
quantum Knizhnik--Zamolodchikov equation of \cite{FR-qKZ}. Note that
this author has not found a prover derivation of this fact in the literature, though it should go along the same lines
as the derivation of the usual $q$KZ equation for matrix elements of VOs of quantized affine algebras
\cite{FR-qKZ,IIJMNT}.

The rational $q$KZ equation is a system of difference equations of the form
\begin{equation}\label{eq:qKZ}
\Psi(z_1,\ldots,z_i+s,\ldots,z_N)
=S_i(z_1,\ldots,z_N)
\Psi(z_1,\ldots,z_i,\ldots,z_N)
,
\qquad i=1,\ldots,N
\end{equation}
where $s$ is equal to the level plus the dual Coxeter number (in units of $\hbar$), so here
$s=(k+1)\hbar$, and the $S_i$ are certain linear operators which
we do no need to describe explicitly.

Instead, one observes that using Eqs.~\eqref{eq:exch} repeatedly (by say moving the argument $z_i+s$ to the
rightmost position in the l.h.s.\ and moving the argument $z_i$ to the leftmost position in the r.h.s.),
one can render \eqref{eq:qKZ}
equivalent to another equation which involves cyclic shift of the parameters $z_i$.
More specifically, one finds the following cyclicity equation:
\begin{equation}\label{eq:cycl}
\Psi^{m_2,\ldots,m_N,m_1}(z_2,\ldots,z_N,z_1+(k+1)\hbar)
=\rho\, \Psi^{m_1,\ldots,m_N}(z_1,\ldots,z_N)
\end{equation}
where the rotation operator $\rho$ is a constant operator (i.e., with no dependence
on the $z_i$) which implements cyclic rotation of factors of the tensor product in 
$L^\ast_{\omega_{m_1}}\otimes \cdots\otimes L^\ast_{\omega_{m_N}}$.
In the case $\lambda=0$, it has an explicit combinatorial definition: 
$\rho_\alpha^\beta=\varepsilon^{m_1}\delta_\alpha^{\rho(\beta)}$,
with $\varepsilon=(-1)^{\frac{M}{k}-1}$, and
$\rho(\beta)$ is obtained from $\beta$ by tableau {\em promotion}\ (evacuation of the $m_1$ letters ``$1$'' followed
by cyclic shift of labels).

Note a technical difference between 
the system (\ref{eq:exch},\ref{eq:cycl}) 
on the one hand and 
\eqref{eq:qKZ} 
on the other hand: the former mixes the various $\Psi^{m_{\sigma(1)},\ldots,m_{\sigma(N)}}$, $\sigma\in\mathcal{S}_N$,
whereas the latter involves a 
single $\Psi^{m_1,\ldots,m_N}$.

\section{Conclusion}
In this short paper,
we have expressed multidegrees of irreducible components of tensor product quiver varieties
in terms of correlators of appropriate Yangian Vertex Operators. By means of bosonization this provides
explicit formulae for them; also, various properties such as conformal block conditions
(see \cite{artic56}) can be easily derived. In a special case, we have shown that they 
are in fact (polynomial) solutions of the $q$KZ equation, though a full geometric proof is lacking. Indeed,
only the exchange relation (the ``nonaffine'' part of $q$KZ, the whole of $q$KZ being related
to an affine Weyl group action) was derived geometrically.

There are many possible generalizations of this work. In particular two underlying Dynkin diagrams can be
generalized away from type A:
\begin{itemize}
\item  The quiver itself, of course, can be chosen of type ADE, leading to
solutions of the Yang--Baxter related to the corresponding double of Yangian;
we also expect the non-ADE solutions (non simply laced or twisted) to be obtainable by appropriate folding.
And one can even go beyond finite Dynkin diagrams, as in the AGT conjecture \cite[part II]{MO-qg}.
\item The gauge group can be chosen of type other than $GL$. In particular type C (symplectic) gauge groups
should lead to systems with boundaries, i.e., solutions of the reflection equation
(see the related work \cite{artic62}), and to solutions of the qKZ equation 
in other types (in the sense of \cite{Che-qKZ}).
However no VO interpretation is known in type other than A.
\end{itemize}

Finally, this work should be extended to K-theory and possibly to elliptic cohomology. Also, it would
be interesting to reinterpret it in the contect of the gauge theory/integrable systems correspondence.

\bigskip
\appendix
\section{A simple example}\label{app:ex}
Consider the quiver of type $A_3$:
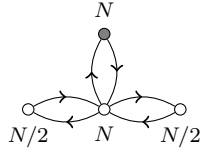
\tikzset{arrow/.style={postaction={decorate,thick,decoration={markings,mark = at position #1 with {\arrow{>}}}}},arrow/.default=0.5}
\tikzset{mynode/.style={circle,draw,inner sep=0mm,minimum size=1.5mm,fill=#1},mynode/.default={white}}
\tikzset{every picture/.style={baseline=-3pt}}
\begin{center}
\begin{tikzpicture}
\node[mynode,label={below:$\ss N/2$}] (a1) at (-1,0) {};
\node[mynode,label={below:$\ss N$}] (a2) at (0,0) {};
\node[mynode,label={below:$\ss N/2$}] (a3) at (1,0) {}; 
\node[mynode=gray,label={above:$\ss N$}] (b) at (0,1) {};
\draw[bend left,arrow] (a1) to (a2);
\draw[bend left,arrow] (a2) to (a1);
\draw[bend left,arrow] (a2) to (a3);
\draw[bend left,arrow] (a3) to (a2);
\draw[bend left,arrow] (a2) to (b);
\draw[bend left,arrow] (b) to (a2);
\end{tikzpicture}
\end{center}
where $N$ is an even integer, the white dots are the $V_a$, the gray dot represents
$W_2$ (all other $W_a$ are zero), and the arrows represent linear maps.

We find
\begin{align*}
\mu&=N \omega_2,  &m&=(\underbrace{2,\ldots,2}_{N})
\\
\lambda&=0,      &\ell&=(\underbrace{4,\ldots,4}_{N/2})
\end{align*}

We claim that the ring of $G_v$-invariant functions has generators
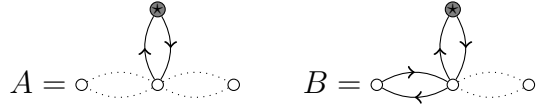
\[
A=
\begin{tikzpicture}
\node[mynode] (a1) at (-1,0) {};
\node[mynode] (a2) at (0,0) {};
\node[mynode] (a3) at (1,0) {}; 
\node[mynode=gray] (b) at (0,1) {$\scriptstyle\star$};
\draw[bend left,arrow] (b) to (a2);
\draw[bend left,arrow] (a2) to (b);
\draw[bend left,dotted] (a1) to (a2);
\draw[bend left,dotted] (a2) to (a1);
\draw[bend left,dotted] (a2) to (a3);
\draw[bend left,dotted] (a3) to (a2);
\end{tikzpicture}
\qquad
B=
\begin{tikzpicture}
\node[mynode] (a1) at (-1,0) {};
\node[mynode] (a2) at (0,0) {};
\node[mynode] (a3) at (1,0) {}; 
\node[mynode=gray] (b) at (0,1) {$\scriptstyle\star$};
\draw[bend left,arrow] (b) to  (a2);
\draw[bend left,arrow] (a2) to (b);
\draw[bend left,arrow] (a1) to  (a2);
\draw[bend left,arrow] (a2) to  (a1);
\draw[bend left,dotted] (a2) to (a3);
\draw[bend left,dotted] (a3) to (a2);
\end{tikzpicture}
\]
(where the picture represents the product of operators along the paths, starting and ending at the marked vertex)
and relations
\[
B(A^2+B)=(A^2+B)B
=A^3+AB+BA=0
\]
which are implied by 
the moment maps conditions
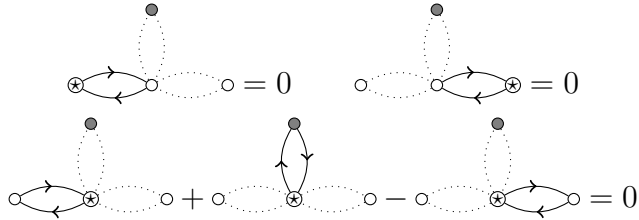
\begin{gather*}
\begin{tikzpicture}
\node[mynode] (a1) at (-1,0) {$\scriptstyle\star$}; 
\node[mynode] (a2) at (0,0) {};
\node[mynode] (a3) at (1,0) {}; 
\node[mynode=gray] (b) at (0,1) {};
\draw[bend left,arrow] (a1) to (a2);
\draw[bend left,arrow] (a2) to (a1);
\draw[bend left,dotted] (a2) to (a3);
\draw[bend left,dotted] (a3) to (a2);
\draw[bend left,dotted] (a2) to (b);
\draw[bend left,dotted] (b) to (a2);
\end{tikzpicture}
=0
\qquad
\begin{tikzpicture}
\node[mynode] (a1) at (-1,0) {};
\node[mynode] (a2) at (0,0) {};
\node[mynode] (a3) at (1,0) {$\scriptstyle\star$}; 
\node[mynode=gray] (b) at (0,1) {};
\draw[bend left,arrow] (a2) to (a3);
\draw[bend left,arrow] (a3) to (a2);
\draw[bend left,dotted] (a2) to (b);
\draw[bend left,dotted] (b) to (a2);
\draw[bend left,dotted] (a2) to (a1);
\draw[bend left,dotted] (a1) to (a2);
\end{tikzpicture}
=0
\\
\begin{tikzpicture}
\node[mynode] (a1) at (-1,0) {};
\node[mynode] (a2) at (0,0) {$\scriptstyle\star$};
\node[mynode] (a3) at (1,0) {}; 
\node[mynode=gray] (b) at (0,1) {};
\draw[bend left,arrow] (a1) to (a2);
\draw[bend left,arrow] (a2) to (a1);
\draw[bend left,dotted] (a2) to (a3);
\draw[bend left,dotted] (a3) to (a2);
\draw[bend left,dotted] (a2) to (b);
\draw[bend left,dotted] (b) to (a2);
\end{tikzpicture}
+
\begin{tikzpicture}
\node[mynode] (a1) at (-1,0) {};
\node[mynode] (a2) at (0,0) {$\scriptstyle\star$};
\node[mynode] (a3) at (1,0) {}; 
\node[mynode=gray] (b) at (0,1) {};
\draw[bend left,arrow] (a2) to (b);
\draw[bend left,arrow] (b) to (a2);
\draw[bend left,dotted] (a2) to (a1);
\draw[bend left,dotted] (a1) to (a2);
\draw[bend left,dotted] (a2) to (a3);
\draw[bend left,dotted] (a3) to (a2);
\end{tikzpicture}
-
\begin{tikzpicture}
\node[mynode] (a1) at (-1,0) {};
\node[mynode] (a2) at (0,0) {$\scriptstyle\star$};
\node[mynode] (a3) at (1,0) {}; 
\node[mynode=gray] (b) at (0,1) {};
\draw[bend left,arrow] (a2) to (a3);
\draw[bend left,arrow] (a3) to (a2);
\draw[bend left,dotted] (a2) to (b);
\draw[bend left,dotted] (b) to (a2);
\draw[bend left,dotted] (a2) to (a1);
\draw[bend left,dotted] (a1) to (a2);
\end{tikzpicture}
=0
\end{gather*}
For example, $A^2+B=\begin{tikzpicture}[scale=0.7]
\node[mynode] (a1) at (-1,0) {};
\node[mynode] (a2) at (0,0) {};
\node[mynode] (a3) at (1,0) {}; 
\node[mynode=gray] (b) at (0,1) {$\sss\star$};
\draw[bend left,arrow] (a2) to (a3);
\draw[bend left,arrow] (a3) to (a2);
\draw[bend left,arrow] (a2) to (b);
\draw[bend left,arrow] (b) to (a2);
\draw[bend left,dotted] (a2) to (a1);
\draw[bend left,dotted] (a1) to (a2);
\end{tikzpicture}
$ and then
\[
B(A^2+B)
=
\begin{tikzpicture}[scale=0.7]
\node[mynode] (a1) at (-1,0) {};
\node[mynode] (a2) at (0,0) {};
\node[mynode] (a3) at (1,0) {}; 
\node[mynode=gray] (b) at (0,1) {$\sss\star$};
\draw[bend left,arrow] (a2) to node[above right=-0.8mm] {$\scriptscriptstyle 6$} (a3);
\draw[bend left,arrow] (a3) to node[below] {$\scriptscriptstyle 7$} (a2);
\draw[bend left,arrow] (a2) to node[above left=-0.8mm] {$\scriptscriptstyle 4,8$} (b);
\draw[bend left,arrow] (b) to node[above right=-0.8mm] {$\scriptscriptstyle 1,5$} (a2);
\draw[bend left,arrow] (a2) to node[below] {$\scriptscriptstyle 2$} (a1);
\draw[bend left,arrow] (a1) to node[above left=-0.6mm] {$\scriptscriptstyle 3$} (a2);
\end{tikzpicture}
=-
\begin{tikzpicture}[scale=0.7]
\node[mynode] (a1) at (-1,0) {};
\node[mynode] (a2) at (0,0) {};
\node[mynode] (a3) at (1,0) {}; 
\node[mynode=gray] (b) at (0,1) {$\sss\star$};
\draw[bend left,arrow] (a2) to node[above right=-0.8mm] {$\scriptscriptstyle 6$} (a3);
\draw[bend left,arrow] (a3) to node[below] {$\scriptscriptstyle 7$} (a2);
\draw[bend left,arrow] (a2) to node[above left=-0.6mm] {$\scriptscriptstyle 8$} (b);
\draw[bend left,arrow] (b) to node[above right=-0.6mm] {$\scriptscriptstyle 1$} (a2);
\draw[bend left,arrow] (a2) to node[below] {$\scriptscriptstyle 2,4$} (a1);
\draw[bend left,arrow] (a1) to node[above left=-0.8mm] {$\scriptscriptstyle 3,5$} (a2);
\end{tikzpicture}
+
\begin{tikzpicture}[scale=0.7]
\node[mynode] (a1) at (-1,0) {};
\node[mynode] (a2) at (0,0) {};
\node[mynode] (a3) at (1,0) {}; 
\node[mynode=gray] (b) at (0,1) {$\sss\star$};
\draw[bend left,arrow] (a2) to node[above right=-1mm] {$\scriptscriptstyle 4,6$} (a3);
\draw[bend left,arrow] (a3) to node[below] {$\scriptscriptstyle 5,7$} (a2);
\draw[bend left,arrow] (a2) to node[above left=-0.6mm] {$\scriptscriptstyle 8$} (b);
\draw[bend left,arrow] (b) to node[above right=-0.6mm] {$\scriptscriptstyle 1$} (a2);
\draw[bend left,arrow] (a2) to node[below] {$\scriptscriptstyle 2$} (a1);
\draw[bend left,arrow] (a1) to node[above left=-0.6mm] {$\scriptscriptstyle 3$} (a2);
\end{tikzpicture}
=0
\]
where the labels indicate the ordering of steps of paths.

For $\lambda=0$, we expect these to be all the relations
(for general $\lambda$, i.e., lower values of the $v$'s, we expect other relations
restricting the ranks of $A$ and $B$).

In order to check this statement, we use the isomorphism with the Mirkovic--Vybornov slice. Here $M=2N$, the nilpotent orbit closure is $\overline{\mathcal O_\ell}=\{X\ M\times M,\ X^4 = 0 \}$, and the slice is
\[
\mathcal T_\mu=\left\{\,
\begin{tabular}{|cc|cc|c|}
\hline
\vrule height11pt depth0pt width0pt 
$0$ & $1$ & $0$ & $0$ & $\cdots$  \\
$\star$ & $\star$ & $\star$ & $\star$ & $\cdots$  \\
\hline
\vrule height11pt depth0pt width0pt 
$0$ & $0$ & $0$ & $1$ & $\cdots$  \\
$\star$ & $\star$ & $\star$ & $\star$ & $\cdots$ \\
\hline
\vrule height11pt depth0pt width0pt 
$\vdots$ & $\vdots$ & $\vdots$ & $\vdots$ & $\ddots$ \\
\hline
\end{tabular}
\,\right\}
\cong
\left\{
\kbordermatrix{&N&N\\ N&0&1\\N&B&A\\}
\right\}
\]
Writing $\begin{pmatrix}
0 & 1 \\
B & A
\end{pmatrix}^4=0$ results in the same equations as before.

In conclusion,
\[
\mathfrak M_0 = \{ A,B\ N\times N\ \text{matrices}:
B(A^2+B)=(A^2+B)B
=A^3+AB+BA=0\}
\]

Similarly, one has
\begin{multline*}
\mathfrak Z_0 = \{ A,B\ N\times N\ \text{strict upper triangular matrices}:
\\
B(A^2+B)=(A^2+B)B
=A^3+AB+BA=0\}
\end{multline*}

$T\cong (\C^\times)^N$ acts by simultaneous conjugation of $A$ and $B$, whereas the extra circle
$\C^\times$ acts by $A\mapsto t^2 A$, $B\mapsto t^4 B$.

Let us now further specialize to $N=4$. Using for example Macaulay2 \cite{M2},
we find that $\mathfrak Z_0$ has three components (which matches
the dimension of $(L_{\omega_2}^{\otimes 4})^{SL(4)}$):
\begin{align*}
\Z_{0,\tinyboxes\tableau{1&2\\1&2\\3&4\\3&4\\}}&=\{ A,B: A_{1,2}=A_{3,4}=B_{1,2}=B_{3,4}=0\}
\\
\Z_{0,\tinyboxes\tableau{1&2\\1&3\\2&4\\3&4\\}}&=\{ A,B: B_{1,2}=B_{2,3}=B_{3,4}=(A^3+AB+BA)_{1,4}=0\}
\\
\Z_{0,\tinyboxes\tableau{1&3\\1&3\\2&4\\2&4\\}}&=\{ A,B: A_{2,3}=B_{2,3}=(AB+BA)_{1,4}=(B^2)_{1,4}=0\}
\end{align*}
with multidegrees
\begin{align*}
\Psi_{\tinyboxes\tableau{1&2\\1&2\\3&4\\3&4\\}}&=(\hbar+z_1-z_2)(\hbar+z_3-z_4)(2\hbar+z_1-z_2)(2\hbar+z_3-z_4)
\\
\Psi_{\tinyboxes\tableau{1&2\\1&3\\2&4\\3&4\\}}&=(2\hbar+z_1-z_2)(2\hbar+z_2-z_3)(2\hbar+z_3-z_4)(3\hbar+z_1-z_4)
\\
\Psi_{\tinyboxes\tableau{1&3\\1&3\\2&4\\2&4\\}}&=(\hbar+z_2-z_3)(2\hbar+z_2-z_3)(3\hbar+z_1-z_4)(4\hbar+z_1-z_4)
\end{align*}

It is not hard to show that the equations
\[
\tau_i\Psi = \check R_i(z_i-z_{i+1})\Psi,\qquad i=1,\ldots,3
\]
define uniquely the $R$-matrices
\begin{align*}
\check R_1(z)=
\check R_3(z)&=
\left(
\begin{array}{ccc}
 \frac{(\hbar-z) (2 \hbar-z)}{(\hbar+z) (2 \hbar+z)} & 0 & 0 \\
 \frac{z (2 \hbar-z)}{(\hbar+z) (2 \hbar+z)} & \frac{2 \hbar-z}{2 \hbar+z} & 0 \\
 \frac{-z(\hbar-z)}{(\hbar+z) (2 \hbar+z)} & \frac{2 z}{2 \hbar+z} & 1 \\
\end{array}
\right)
\\
\check R_2(z)&=\left(
\begin{array}{ccc}
 1 & \frac{2 z}{2 \hbar+z} & \frac{-z(\hbar-z)}{(\hbar+z) (2 \hbar+z)} \\
 0 & \frac{2 \hbar-z}{2 \hbar+z} & \frac{z (2 \hbar-z)}{(\hbar+z) (2 \hbar+z)} \\
 0 & 0 & \frac{(\hbar-z) (2 \hbar-z)}{(\hbar+z) (2 \hbar+z)} \\
\end{array}
\right)
\end{align*}
and that these satisfy Eqs.~(\ref{eq:ybe}--\ref{eq:comm}).

Furthermore, $\Psi$ satisfies the remarkable cyclicity relation (Eq.~\eqref{eq:cycl})
\[
\Psi(z_2,z_3,z_4,z_1+5\hbar)=\rho\,\Psi(z_1,z_2,z_3,z_4)
\]
where $\rho=\begin{pmatrix}0&0&1 \\ 0&1&0 \\ 1&0&0
\end{pmatrix}$. (the fact that $\rho^2=1$ explains that $\check R_1(z)=\check R_3(z)$.)

Finally, $\Psi$ satisfies the wheel condition:
\[
\Psi(z,z+2\hbar,z+4\hbar,\cdot)
=
\Psi(z,z+2\hbar,\cdot,z+4\hbar)
=
\Psi(z,\cdot,z+2\hbar,z+4\hbar)
=
\Psi(\cdot,z,z+2\hbar,z+4\hbar)
=
0
\]

\rem{a comment on connection to Brauer loop scheme?}

In order to switch to the stable basis, one must deform $\M_0$ and $\Z_0$
by taking the fiber at generic values of the moment map, or equivalently
by deforming the nilpotent orbit closure $\{X^4=0\}$ to the regular orbit
$\{\prod_{a=1}^4 (X-t_a)=0,\ \det(x-X)=\prod_{a=1}^4(x-t_a)^{\lambda_a}\}$
\cite{articzz}. The equations deform into
\begin{align*}
A^2 B + B^2 - e_1 A B + e_2 B + e_4 &= 0
\\
B A^2 + B^2 - e_1 B A + e_2 B + e_4 &= 0
\\
A^3 + A B + B A - e_1 (A^2+B) +e_2 A - e_3 &= 0
\end{align*}
where the $e_i$ are the elementary symmetric polynomials of the $t_i$.
The deformation $\Z_t$ of $\Z_0$ corresponds to taking $A$ and $B$ to
be upper triangular (not strict upper triangular), the diagonal entries
encoding the stable basis components; in the present case, the
eigenvalues of $A$ are $t_i+t_j$, $1\le i<j\le 4$, corresponding to
the $GL(4)$ weights on $L_{\omega_2}$. 
In order to select the weight space,
one must either keep only the top-dimensional components, or impose
the extra constraint on the characteristic polynomial of $X$ mentioned above.

The coefficients of the change of basis are obtained by taking the flat
limit $t_i\to 0$: each limit $\Z'_{0,\alpha}$ of an irreducible component $\Z'_{t,\alpha}$
of $\Z_t$ is a union of $\Z_{0,\beta}$ with multiplicity $|c^\beta_\alpha|$.
Here we shall not write down all of the $c^\beta_\alpha$ explicitly because already in size $4$,
the zero weight space of $L_{\omega_2}^{\otimes 4}$ has dimension $90$.
Instead we give only one component corresponding to $\alpha=(\{1,2\},\{1,3\},\{2,4\},\{3,4\})$:
\begin{align*}
\Z'_{t,\alpha}=\big\{A,B:\ &A_{i,i}=\sum_{a\in\alpha_i} t_a, B_{i,i}=-\prod_{a\in\alpha_i} t_a, 
\\
&B_{2,3}A_{3,4}+(t_3-t_2)B_{2,4}+t_3 A_{2,3}A_{3,4}+t_3(t_3-t_2)A_{2,4}=0
\\
&A_{1,2}B_{2,3}+(t_2-t_3)B_{1,3}+t_2 A_{1,2}A_{2,3}+t_2(t_2-t_3)A_{1,3}=0
\\
&A_{1,2}B_{2,4}+B_{1,3}A_{3,4}+A_{1,2}A_{2,3}A_{3,4}+t_2A_{1,3}A_{3,4}+t_3A_{1,2}A_{2,4}=0
\big\}
\end{align*}
The $t_i\to 0$ limit $\Z'_{0,\alpha}$ is the union of 
$\Z_{0,\tinyboxes\tableau{1&2\\1&2\\3&4\\3&4\\}}$ and
$\Z_{0,\tinyboxes\tableau{1&2\\1&3\\2&4\\3&4\\}}$, each with multiplicity one.

\gdef\MRshorten#1 #2MRend{#1}%
\gdef\MRfirsttwo#1#2{\if#1M%
MR\else MR#1#2\fi}
\def\MRfix#1{\MRshorten\MRfirsttwo#1 MRend}
\renewcommand\MR[1]{\relax\ifhmode\unskip\spacefactor3000 \space\fi
\MRhref{\MRfix{#1}}{{\scriptsize \MRfix{#1}}}}
\renewcommand{\MRhref}[2]{%
\href{http://www.ams.org/mathscinet-getitem?mr=#1}{#2}}
\bibliographystyle{amsplainhyper}
\bibliography{biblio}
\end{document}

%% file: tableau.tex
\newdimen{\cellsize}
\newcommand\medboxes{\setlength{\cellsize}{14pt}\def\boxformat{}}

\newcommand\tinyboxes{\setlength{\cellsize}{6pt}\def\boxformat{\scriptscriptstyle}}
\medboxes
\tikzset{tableaubox/.style={draw=black,thin,sharp corners,solid,minimum size=\cellsize,inner sep=0pt}}
\tikzset{tableau/.style={matrix,name=tab,matrix anchor=tab-1-1.south west,inner sep=1pt,matrix of math nodes,cells={anchor=center,draw=black,thin,solid,arrows=-},nodes={tableaubox,execute at begin node=\boxformat},nodes in empty cells,row sep={\cellsize,between origins},column sep={\cellsize,between origins}}}
\newcommand\missingcell{|[draw=none]|}
\makeatletter
\newcommand\cellextra[1]{#1\expandafter\tikz@lib@matrix@start@cell}
\makeatother
\newcommand\vdotscell{\cellextra{\draw[dotted] (-0.5*\cellsize,-0.5*\cellsize) --++(0,\cellsize);}\missingcell}
\newcommand\vcell{\cellextra{\draw (-0.5*\cellsize,-0.5*\cellsize) --++(0,\cellsize);}\missingcell}
\def\activate#1{\begingroup
  \lccode`\~=`#1%
  \lowercase{\endgroup \let~#1}%
  \catcode`#1=13\relax}
\activate &
\def\tableau#1{\tikz[baseline=0]
\node[tableau]{#1};%
}